\documentclass[12pt]{article}

\usepackage{amsmath}
\usepackage{amsfonts}
\usepackage{amssymb, amsfonts, graphicx, amsmath, nicefrac}
\usepackage{bbm}
\usepackage{xcolor}
\usepackage{stmaryrd}
 \usepackage{setspace}
\newcommand{\ccr}{{\mathfrak r}}
\newcommand{\cm}{{\sf m}}

\newtheorem{definition}{Definition}%[section]

\newcommand{\N}{{\cal N}}

\newcommand{\T}{{\cal T}}

\newcommand{\G}{{\cal G}}

\newcommand{\M}{{\cal M}}
\newcommand{\X}{{\cal X}}

\newcommand{\be}{\begin{equation}}
\newcommand{\ee}{\end{equation}}

\def\O{\Omega}

\def\rank{{\rm rank}}
\def\dim {{\rm dim}}

\def\bbr{{\Bbb{R}}} %real numbers
\newcommand{\yao}[1]{{\color{black}{#1}}}

\newcommand{\rui}[1]{{\color{black}{#1}}}

\usepackage{mathtools}
 \usepackage{nccmath}
 \usepackage{cite}

\title{On characteristic rank for  matrix and tensor completion}
\author{Alexander Shapiro, Yao Xie, and Rui Zhang}
\date{\today}

\begin{document}

\maketitle

\section{Scope}

In this lecture note,  we discuss a fundamental concept, referred to as the {\it characteristic rank}, which suggests a general framework for characterizing the basic properties of various low-dimensional models used in signal processing.  Below, we illustrate this framework using two examples: matrix and three-way tensor completion problems, and consider basic properties include identifiability of a matrix or tensor, given partial  observations. In this note, we consider cases without observation noise to illustrate the principle.

\section{Relevance}

Characteristic rank provides a fundamental tool for determining the ``order'' of low-rank structures, such as the rank of low-rank matrices and rank of three-way tensors. The concept of characteristic rank was introduced in \cite{sxz19}, where it was used to establish necessary and sufficient conditions to determine the ``recoverability'' of the low-rank matrices in \cite{sxz19}.

\vspace{.1in}

The characteristic rank can also be generally applied to determine the ``intrinsic'' degree-of-freedom in other low-rank manifold structures. Such instances include determining the number of hidden nodes in one-layer neural networks and determining the number of sources in blind demixing problems, as shown in \cite{Shapiro2019GoodnessoffitTO}.

\section{Prerequisite}

To better comprehend the concepts discussed in this lecture-notes article, readers are expected to have a good background in linear algebra, multivariate calculus, and basic concepts of measure theory (which we will explain whenever running into them). Suggested references are \cite{rudin1964principles} and \cite{meyer2000matrix}. Below, we review some basic concepts necessary for the notes.

\vspace{.1in}
\noindent
{\bf Manifold of low-rank matrices.}
Consider  the set   of $n_1\times n_2$   matrices of rank $r$, denoted $\M_r$. Note that the rank is no larger than the dimension of the matrix: $r\le \min\{n_1,n_2\}$. It is known that such a set of rank-$r$ matrices $\M_r$
forms a {\it smooth manifold}  in the space $\bbr^{n_1\times n_2}$ and the dimension of the manifold is given by
  \begin{equation}\label{dimrank}
  \dim(\M_r)=r(n_1+n_2-r).
  \end{equation}
A matrix $A\in \M_r$    can be represented in the form
$A=V W^\top,$
where $V$ and $W$ are matrices of the respective order $n_1\times r$  and $n_2\times r$, both of full column rank $r$. Thus, we can view $(V,W)$ as a parametrization of $\M_r$. Note that the number of involved parameters is $r(n_1+n_2)$, which is larger than the dimension of $\M_r$; this is because $V$ and $W$  in  the  above representation are  not unique.

\vspace{.1in}
\noindent
{\bf Three-way tensor.} Another example we will consider is the (three-way) tensor $X\in \bbr^{n_1\times n_2\times n_3}$. It is said that $X$  has rank one if
$X=a\circ b\circ c$  where $a,b,c$ are  vectors of the  respective dimensions $n_1,n_2,n_3$, and $``\circ"$  denotes the vector outer product. That
is, every element of tensor $X$  can be written as the product
$X_{ijk}=a_i b_j c_k$.  The smallest number $r$  such that tensor $X$ can be represented as a sum of $r$  of rank-one tensors is called the rank of $X$. The corresponding
decomposition is often referred to as the (tensor) rank-decomposition or Canonical Polyadic
Decomposition (CPD) \cite{sorber2013optimization,sidiropoulos2017tensor}. \rui{We would like to remark that our method can apply to higher order tensors as well. }

\section{Problem statement}

\subsection{Matrix completion}
% \vspace{.1in}
% {\bf Matrix Completion.}
Let us start by considering the  problem of reconstructing  an $n_1\times n_2$ matrix  of a given rank $r$ while observing  its entries $M_{ij}$, $(i,j) \in \O$,  for an index set $\O\subset \{1,...,n_1\}\times \{1,...,n_2\}$ of cardinality $\cm=|\O|$. This is known as the exact matrix completion problem \cite{candes2010power}, which is now well-studied. The conditions for recovery have been derived assuming entries are {\it missing-at-random} and the performance guarantees are given in a probabilistic sense. Here, we aim to approach the problem from a geometric perspective, which can possibly lead to a deterministic and more intuitive answer.  There are two basic problems associated with this problem, namely the {\it existence} and {\it uniqueness} of the solution. That is,  whether such a matrix does exist and if it exists, whether it is unique. Fundamentally, these questions are related to the identifiability of low-rank matrices, which we define as follows.
\begin{definition}[Local identifiability of low-rank matrix \rui{completion problem}] Let $Y\in \M_r$  be such that $[Y]_{ij}=M_{ij}$,  $(i,j) \in \O$. (Thus, $rank(Y)=r$.) It is said that the matrix completion problem is locally identifiable at $Y$ if there exists a neighborhood $\N\subset \bbr^{n_1\times n_2}$ of $Y$ such that for any $Y'\in \N$ with $[Y']_{ij}=M_{ij}$,  $(i,j) \in \O$, the rank of $Y'$ is different from $r$.
\end{definition}

\subsection{Uniqueness of tensor decomposition}
% \vspace{.1in}
% %
% \noindent{\bf Uniqueness of tensor decomposition.}
Uniqueness is the key question of the tensor rank decomposition. Here, we consider the following tensor decomposition problem: given a three-way tensor $X$, we would like to find the associated matrix factors $A, B, C$  of the respective order
$n_1\times r$, $n_2\times r$  and $n_3\times r$, such that \rui{$X = \llbracket A,B,C\rrbracket$}, meaning that $X=\sum_{i=1}^r  a^i\circ b^i\circ c^i$   with
$a^i,b^i,c^i$ being $i$th columns of the respective matrices $A,B,C$.
Clearly the decomposition \rui{$X = \llbracket A,B,C\rrbracket$}
  is invariant with respect to permutations of the rank one components,
and rescaling of the columns of matrices $A,B,C$  by factors
$\lambda_{1i},\lambda_{2i},\lambda_{3i}$
such that $\lambda_{1i}\lambda_{2i}\lambda_{3i}=1$, $i=1, \ldots ,r$. We first introduce the global and local identifiability of tensor.

\begin{definition}[Global identifiability of tensor \rui{decompositon}]
The decomposition \rui{$X = \llbracket A,B,C\rrbracket$} is (globally) {\em identifiable} of rank $r$  if it is unique, i.e., if \rui{$X = \llbracket A',B',C'\rrbracket$} is another decomposition of tensor  $X$  with matrices  $A',B',C'$  being of the respective order
$n_1\times r'$, $n_2\times r'$, $n_3\times r'$ and $r'= r$, then both decompositions are the same up to the corresponding permutation
and rescaling. It is said that the rank $r$  decomposition is {\em generically} identifiable if for almost every $(A,B,C)\in \bbr^{n_1\times r}\times \bbr^{n_2\times r}\times \bbr^{n_3\times r}$  the corresponding tensor
 $X=\llbracket A,B,C\rrbracket $ is identifiable of rank $r$.
 \end{definition}

% \textcolor{red}{(Yao: Do we need to define global idenfiability here? It seems that the following discussions do not need this concept or do I need anything? )}

\begin{definition}[Local identitifiability of tensor \rui{decomposition}]
\label{def-lociden}
We say that $(A,B,C)\in  \bbr^{n_1\times r}\times   \bbr^{n_2\times r}\times \bbr^{n_3\times r} $ is   {\em  locally  identifiable}
if there is a neighborhood $\N$ of $(A,B,C)$ such that
$(A',B',C')\in \N$ and \rui{$\llbracket A',B',C'\rrbracket =  \llbracket A,B,C\rrbracket$} imply  that $(A',B',C')$ can be obtained from $(A,B,C)$ by the corresponding rescaling.
We say that model $(n_1,n_2,n_3,r)$   is  {\em generically} locally   identifiable if a.e. $(A,B,C)\in  \bbr^{n_1\times r}\times   \bbr^{n_2\times r}\times \bbr^{n_3\times r} $ is  locally  identifiable.
\end{definition}

Like the matrix completion problem,  it is also possible to consider a tensor completion problem: reconstructing a tensor of a given rank when only a subset of the entries is observed. The respective local and global identifiability concepts can be defined similarly.
%Various studies are giving sufficient upper bounds, ensuring that if  $r$ is less than those bounds, then the rank $r$  decomposition is generically identifiable. However, the respective exact upper bounds are not known.

\section{Solutions}

\subsection{Matrix completion}

{\bf Re-parameterization of matrix completion problem.} Let us start with the matrix completion problem using the following parametrization.
Consider the set  $\X$   of $n_1\times n_2$  matrices $X$ such that $[X]_{ij}=0$,  $(i,j) \in \O$ (imagining adding such matrices to solutions and they are still consistent with observations). We can view $\X$ as a linear space of dimension $\dim(\X)=n_1 n_2-\cm$. Then the  matrix completion problem has a solution if and only if  there exist respective matrices $V$ and $W$ of rank $r$  and $X\in \X$ such that
   $[V W^\top +X]_{ij}=M_{ij}$,  $(i,j) \in \O$.
Let $\Theta$ be the set of vectors $\theta$ formed from the components of    $(V,W,X)$. Note that $\Theta$ is a subset of vector space of dimension $r(n_1+n_2)+n_1 n_2-\cm$.

\vspace{.1in}
\noindent {\bf Characteristic rank.}
The  matrix completion  parametrization  can be considered as a mapping assigning matrix $V W^\top +X$ to vector of parameters $\theta=(V,W,X)\in \Theta$. With this mapping, we can define the so-called {\it Jacobian matrix} $\Delta(\theta)$, which is the partial derivatives of $V W^\top +X$
with respect to components of vector $\theta$. Then we associate this mapping with its  {\em characteristic rank}, defined as
 \begin{equation}\label{eq-3}
  \ccr=\max_{\theta\in \Theta} \{\rank (\Delta(\theta))\}.
 \end{equation}
 Note that the characteristic rank $\ccr$ does not depend on order in which the parameters are arranged.

\vspace{.1in}
 The characteristic rank has the following properties: the rank of
  $\Delta(\theta)$ is equal to $\ccr$ for almost every (a.e.) $\theta\in \Theta$. By ``almost every" we mean that the set of such $\theta\in \Theta$ for which  $\rank (\Delta(\theta))\ne \ccr$ is of Lebesgue measure zero. Moreover, the
set $\{\theta\in \Theta: \rank (\Delta(\theta))=\ccr\}$ forms an open subset of $\Theta$. It follows that  the rank of
  $\Delta(\theta)$ is constant, and  equals $\ccr$,  in a neighborhood of a.e. $\theta\in \Theta$. \yao{ This result implies that {\it the characteristic rank is an intrinsic quantity associated with the ``degrees-of-freedom'' of the problem, regardless of the value of the parameters.}}

\vspace{.1in}
\noindent  {\bf Implication of characteristic rank on matrix completion.}
 We can also look at the characteristic rank from the following point of view. Consider the tangent space $\T_{\M_r}(Y)$ to the manifold $\M_r$, at the point $Y=V  W^\top\in \M_r$.  We have that
 \begin{equation}\label{rank-1}
  \rank(\Delta(\theta))=\dim (\T_{\M_r}(Y))+\dim (\X)-\dim\left (\T_{\M_r}(Y)\cap \X\right),
 \end{equation}
The above relation (\ref{rank-1}) can be explained as follows. Generically the image of the considered mapping $V W^\top +X$ forms a smooth manifold in the image space, at least locally.  The tangent space to this manifold, at the considered point, is the sum of the tangent space to $\M_r$ (from the paramerization $VW^\top$) and the linear space $\X$ in the image space. On the other hand this tangent space is generated by columns of the Jacobian matrix $\Delta(\theta)$ (or in other words by the differential of the mapping) and its dimension is equal to the rank of  $\Delta(\theta)$.  Then the right hand side of (3) is the usual formula for dimension of the sum of two linear spaces $\T_{M_r}(Y)$ and $\X$.
Hence, from (\ref{rank-1}) and the definition of the characteristic rank (\ref{eq-3}), we have that
 \begin{equation}\label{rank-2}
  \ccr=\dim (\T_{\M_r}(Y))+\dim (\X)-\inf_{Y\in \M_r}\left\{\dim\big (\T_{\M_r}(Y)\cap \X\big)\right\}.
 \end{equation}
%
%\begin{equation}\label{eq-4}
%  \ccr\le \dim (\T_{\M_r}(Y))+\dim (\X)=r(n_1+n_2-r) + n_1 n_2-\cm.
% \end{equation}
By classical  Sard's theorem \cite{sards1942measure}, we have  that the image of the set  $\Theta$ by  the mapping $\theta\mapsto V W^\top +X$,  has Lebesgue measure zero if and only if $\ccr< n_1  n_2$. That is, if $\ccr< n_1  n_2$, then {\em generically} the problem of  reconstructing   matrix of rank $r$ by observing its entries  $M_{ij}$, $(i,j) \in \O$, is unsolvable. By  ``generically" we mean that the set of rank-$r$ solutions with components matching $M_{ij}$, $(i,j) \in \O$, has a Lebesgue measure zero in the corresponding vector space of dimension $\cm$.

\vspace{.1in}
In other words, if the characteristic rank is smaller than the dimension $n_1  n_2$ of the image space, then any solution of rank $r$ is unstable: this means that arbitrarily small changes of the data values $M_{ij}$ make rank $r$ solution unattainable.  Note that the characteristic rank is a function of the index set $\O$ and does not depend on the observed values $M_{ij}$. In particular,
because of \eqref{rank-2} we have that $\ccr< n_1n_2$ if
$
 \cm > r(n_1+n_2-r) .
$
For example, if $n_1 = n_2 = 10$, $r = 3$, then we have $\ccr < 100$ if $\cm > 3\times (10+10-3) = 51$. Since the characteristic rank is the dimension of the image of the mapping, if it is smaller than the dimension \rui{$n_1 n_2$} of the image space, then it is ``thin", i.e. of measure zero in the image space. %Formula (4) follows directly from equations  (3) and (2).
%\yao{(Question: what's the implication of this? When the number of observation is larger than the lower bound (which depends on $r$ and the dimension of the matrix, then we cannot uniquely fill out the entries of a matrix such that it matches the observation, and the rank of the matrix is $r$?) For example, if $n_1 = n_2 = 10$, $r = 3$, then we cannot have a unique solution if $\cm > 3\times (10+10-3) = 51$, but when $\cm < 51$, we can find a unique solution? Thanks.)}

\vspace{.1in}
\noindent{\bf Well-posedness condition.}
By the above discussion we have that if
 \begin{equation}\label{rank-3}
 \T_{\M_r}(Y)\cap \X=\{0\}
 \end{equation}
 at least for one point $Y\in \M_r$, then
 \begin{equation}\label{rank-4}
 \ccr=\dim (\T_{\M_r}(Y))+\dim (\X).
 \end{equation}
 Conversely if  \eqref{rank-4} holds, then condition \eqref{rank-3} is satisfied for all $Y\in \M_r$ except for a  set of measure zero in  $\M_r$. Condition   \eqref{rank-3} implies {\it local identifiability} at $Y$.

 \begin{itemize}
   \item
   Generically the matrix completion problem is locally identifiable if and only if condition  \eqref{rank-4} holds, which is referred to as the {\it well-posedness condition} in \cite{sxz19}.
 \end{itemize}

\begin{figure}[h!]
    \centering
   \includegraphics[width = 0.4\linewidth]{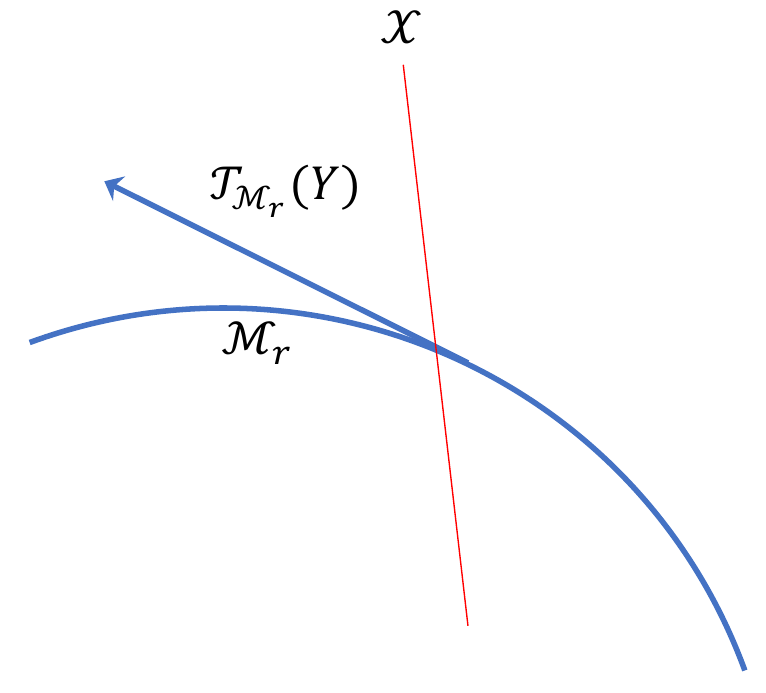}
    \caption{Illustration of well-posedness condition for the matrix completion problem.}
 \label{wellpose_graph}
 \end{figure}

 Figure \ref{wellpose_graph} illustrates the above point. Generically  the intersection of $\T_{M_r}(Y)$ and $\X$ gives the tangent space to the intersection of  $M_r$  and $\X$.   When the intersection of $\T_{M_r}(Y)$ and $\X$ is $\{0\}$ we have well posedness and local uniqueness.

  \vspace{.1in}
\noindent  {\bf Simple example.}
  Here we illustrate the characteristic rank using a simple example of 2-by-2 rank-one matrix  $M = vw^\top$, with partial observations at $\Omega = \{(1,1), (2,2)\}$. Then
\[X = \begin{bmatrix} 0 & x_{12} \\ x_{21} & 0
\end{bmatrix},\]
$\theta = (v_1, v_2, w_1, w_2, x_{12}, x_{21})$. We have,
\begin{equation*}
    \Delta(\theta) = \frac{\partial (vw^\top + X)}{\partial\theta} =
\begin{bmatrix}
     w_1& 0   &v_1&0  &0&0  \\
     w_2& 0   &0  &v_1&1&0  \\
     0  & w_1 &v_2&0  &0&1  \\
     0  & w_2 &0  &v_2&0&0
\end{bmatrix}.
\end{equation*}
It can be verified that $\rank (\Delta(\theta)) = 4$ for a.e. $\theta \in \Theta$; thus, $\ccr = 4$. Consider possible rank-one solution to this problem. The tangent space of the rank-one manifold $\dim (\T(\M_r)) = 2+2 - 1 = 3$, and $\dim (\X) = 2$; $\ccr < \dim (\T(\M_r))+\dim (\X)$ and the well-posedness condition (\ref{rank-4}) is not satisfied. Indeed, the rank-one solution to this problem is not unique: it can be any $x_{12} x_{21} = c$ where c is the product of the observed diagonal elements.

% On the other hand, if $\Omega = \{(1,1), (1,2)\}$, \[X = \begin{bmatrix} 0 & 0 \\ x_{21} & x_{22}
% \end{bmatrix},\] $\theta = (v_1, v_2, w_1, w_2, x_{21}, x_{22})$. We have,
% \begin{equation*}
%     \Delta(\theta) = \frac{\partial (vw^\top + X)}{\partial\theta} = \left(
% \begin{array}{cccccc}
%      w_1& 0   &v_1&0  & 0&0  \\
%      w_2& 0   &0  &v_1& 0&0  \\
%      0  & w_1 &v_2&0  & 1&0  \\
%      0  & w_2 &0  &v_2& 0&1
% \end{array}.\right)
% \end{equation*}
% It can be verified that $\rank (\Delta(\theta)) = 4$ regardless of the value of $\theta$.

On the other hand, if $\Omega = \{(1,1), (1,2), (2,1)\}$, \[X = \begin{bmatrix} 0 & 0 \\ 0 & x_{22}
\end{bmatrix},\] $\theta = (v_1, v_2, w_1, w_2, x_{22})$. We have,
\begin{equation*}
    \Delta(\theta) = \frac{\partial (vw^\top + X)}{\partial\theta} =
    \begin{bmatrix}
     w_1& 0   &v_1&0  &0   \\
     w_2& 0   &0  &v_1&0  \\
     0  & w_1 &v_2&0  &0  \\
     0  & w_2 &0  &v_2&1
\end{bmatrix}.
\end{equation*}
It can be verified that $\rank (\Delta(\theta)) = 4$ for a.e. $\theta \in \Theta$, and thus $\ccr = 4$. The rank of the tangent space is $2+2 -1=3$, the dimension of $\X$ is 1. Thus, $ \ccr=\dim (\T(\M_r))+\dim (\X)$ and the well-posedness condition (\ref{rank-4}) is satisfied. Indeed, the solution to this matrix completion problem is unique.

\vspace{.1in}
\noindent
{\bf Checking conditions.} Although the above simple example is easy to check, to evaluate the characteristic rank in a closed-form is not always easy for larger instances. Nevertheless, the rank of the Jacobian matrix can be computed numerically, and hence condition  \eqref{rank-4} can be verified for a considered index set $\O$ and rank $r$.  Clearly, local identifiability is a necessary condition for global identifiability (i.e., for global uniqueness of the solution). Assuming that all observed entries are different from zero,    necessary and sufficient conditions for global identifiability are known when  $r=1$. Those conditions are the same for local identifiability \yao{(see \cite{sxz19} for more details)}. To give necessary and sufficient conditions for global identifiability for general $r$ and $\O$ could be too difficult and out of reach. On the other hand,  the simple dimensionality condition \eqref{rank-4} gives a verifiable condition at least for local identifiability.

\subsection{Tensor decomposition}

\vspace{.1in}
\noindent{\bf Invoking characteristic rank on three-way tensor.} Here we briefly discuss the local identifiability for tensor decomposition. For three-way tensor recovery, we can consider the mapping
 \begin{equation}\label{tens-1}
\rui{\G_r}:(A,B,C)\mapsto  \rui{\llbracket A,B,C\rrbracket}.
 \end{equation}
Similar to \eqref{eq-3},  the characteristic rank $\ccr$  of the above mapping is given by the maximal rank of its Jacobian matrix, and it has generic properties similar to the ones discussed for the matrix completion problem.
 Note that $\ccr$ is always less than or equal to $r(n_1+n_2+n_3-2)$. This follows by counting the number of elements in  $(A,B,C)$
and making correction for the scaling factors.
\begin{itemize}
  \item
The model $(n_1,n_2,n_3,r)$   is   generically  locally   identifiable if and only if the following condition for the  characteristic rank holds
\begin{equation}\label{tens-2}
 \ccr=r(n_1+n_2+n_3-2).
\end{equation}
\end{itemize}
The above condition \eqref{tens-2}  is necessary for the generic  global identification, and can be verified numerically by computing rank of the Jacobian matrix of the mapping $\G_r$.

\vspace{.1in}
Let us note further that in a similar spirit, it is also possible to give conditions for local identifiability of the tensor completion problem when only a set of observed values of the tensor components are available \yao{(i.e., tensor completion problems)}. To do so, we need to set up appropriate mapping and study the associated characteristic rank.

We can refer to \cite[section 3.2]{kolda2009tensor}, and references there in,  for a discussion of uniqueness (identifiability) of tensor rank decompositions. For the tensor completion problem, local identifiability does not imply the respective global identifiability even in the rank one case (e.g., \cite{Mohit}).

\section{Computational example}

\begin{figure}[h!]
  \centering
    \includegraphics[width = 0.6\linewidth]{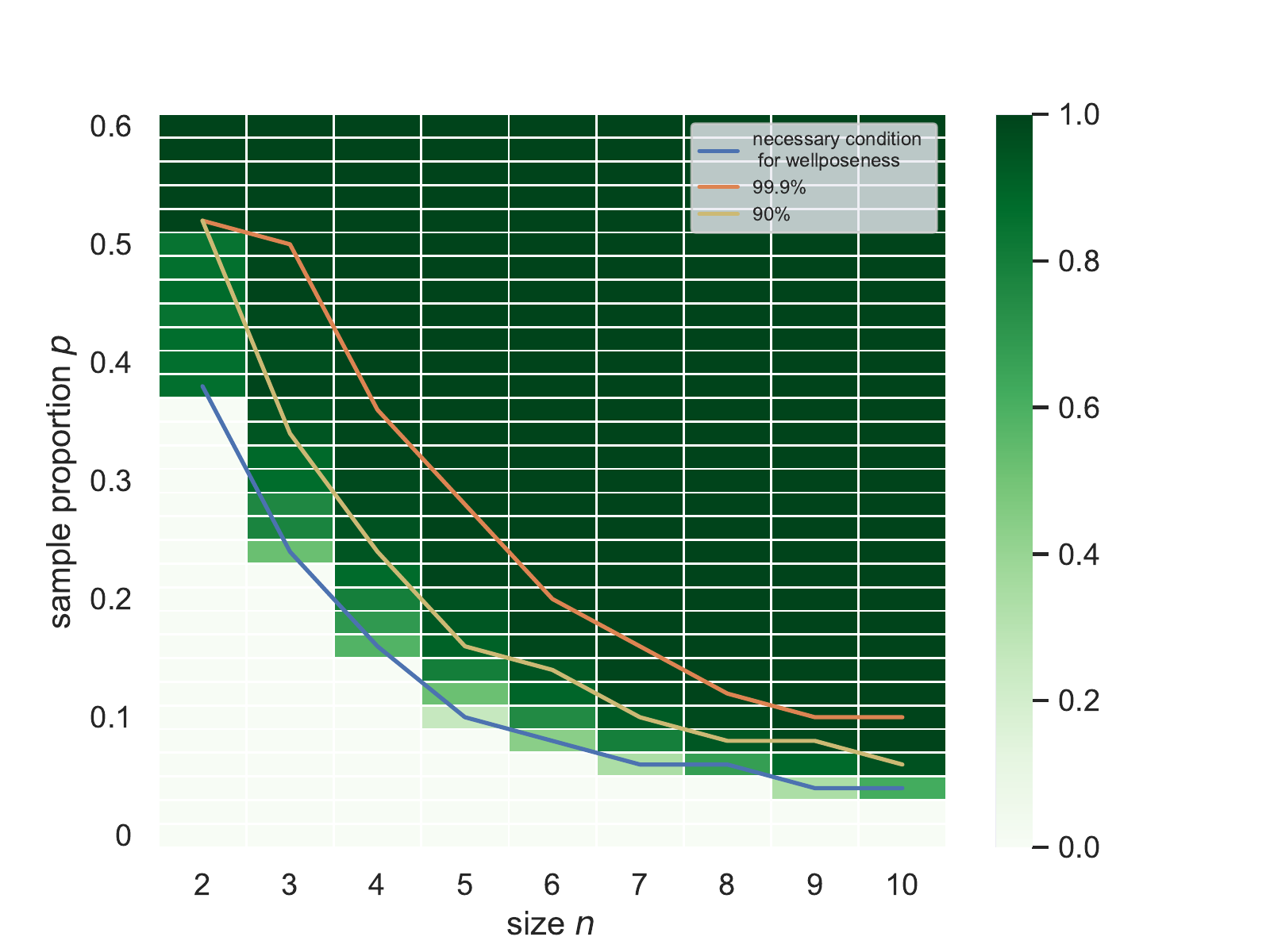}
    \caption{Example of the recovering a three-way tensor with missing data: the probability of well-posedness being satisfied versus theoretical prediction. The blueline corresponds to $p = (3n-2)/n^3$; the yellow and the orange lines correspond to the sampling proportion that the well-posedness condition is satisfied with probability $90\%$ and $99.9\%$ empirically.}
  \label{wellpose}
\end{figure}

\rui{Here we present a numerical example to illustrate how to use the characteristic rank to study a three-way tensor's completion problem. Consider the case where the tensor entries are randomly sampled. Assume the size of each dimension of the tensor is $n$, and thus the size of the tensor is $\bbr^{n\times n\times n}$. The proportion of the observed  entries is $p$,  and the total number of observed entries is  $\cm = \lceil pn^3 \rceil$, where $\lceil x\rceil$ is the ceiling function for rounding up to the nearest integer. For each $p$, we randomly choose $m$ observations from the tensor. For the reported experiments, we used  $n=2,\dots, 10$ and $p = 0.02, 0.04, \ldots, 0.6$. As previously mentioned, the necessary condition for the well-posedness is that $\cm \ge 3n-2$ (see (\ref{tens-2})). This requires, approximately, $p \ge (3n-2)/n^3$. Figure \ref{wellpose} shows the probability that the well-posedness is satisfied for rank-one tensors under different tensor sizes and sampling proportions. Note that the empirical results match well with the theoretical prediction. Moreover, it can be observed that as the tensor size becomes large, the well-posedness condition is satisfied with a small sampling proportion.}
% Here we present a numerical example to illustrate how to use the characteristic rank to study a three-way tensor's completion problem. Consider the case where the tensor entries are randomly missing with probability $p$. In the set of experiments, the size of tensor is $\bbr^{n\times n\times n}$, and each element is sampled with probability $p$, where $n=2,\dots, 10$ and $p=0.05, 0.1,\dots, 1$. Since $\mathbb E[\rui{\cm}] = n^3 p$, we have $p\approx \rui{\cm}/n^3$. As previously mentioned, a necessary condition of well-posedness, following (\ref{tens-2}), is that $\cm > 3n-2$. This requires, roughly, $p > (3n-2)/n^3$. Figure \ref{wellpose} shows the probability that the well-posedness is satisfied for rank-one tensor under different tensor size and sampling probability for randomly missing entries. Note that the empirical results match well with the theoretical prediction. Moreover, it can be observed that as the tensor size becomes large, the well-posedness condition is satisfied with small sampling probability.

%   \begin{figure}[h!]
%   \centering
%     \includegraphics[width = 0.6\linewidth]{10_0.05_300_transition.pdf}
%     \caption{Example of the recovering a three-way tensor with missing data: the probability of well-posedness being satisfied versus theoretical prediction. The blueline corresponds to $p = (3n-2)/n^3$; yellow and orange lines correspond to the sampling probability that well-posedness condition is satisfied with probability $90\%$ and $99.9\%$ empirically.}
%   \label{wellpose}
%  \end{figure}

\section{What we have learned}

In this note, we explained how to use a fundamental concept, namely the {\it characteristic rank}, to answer essential questions such as {\it identifiability} when given observations of a low-rank structure (e.g., low-rank matrices and low-rank three-way tensors). The framework involves a few steps. We first find the map that associates the truth to the observations, then study the Jacobian matrix of the map to find the characteristic rank, and compare the characteristic rank with respective conditions that are problem-specific (such as well-posedness condition). Once the concepts are understood, the analysis usually involves only basic multi-variate calculus. The benefit is that the tool can generally be applicable to study other problems with low-rank structures. In this note, we have considered cases without observation noise to illustrate the principle. When there are additive Gaussian noises, statistical goodness-of-fit tests can be developed based on the framework  \cite{Shapiro2019GoodnessoffitTO}.

\section*{Author Biography}

Alexander Shapiro is  A. Russell Chandler III Chair and Professor   at Georgia Institute of Technology in the H. Milton Stewart School of Industrial and Systems Engineering. He has published more than 140 research articles in peer review journals and is a coauthor of several
books. He 
served on editorial board of a number of professional journals. He was the Editor-In-Chief of the Mathematical Programming, Series A, journal. In 2013 he was a recipient of Khachiyan Prize for Life-time Accomplishments in Optimization, awarded by the INFORMS Optimization Society, and in 2018 he was a recipient of the Dantzig Prize awarded by the Mathematical Optimization Society and Society for Industrial and Applied Mathematics. In 2020 he was elected to the National Academy of Engineering.

Yao Xie is an Associate Professor and Harold R. and Mary Anne Nash Early Career Professor at Georgia Institute of Technology in the H. Milton Stewart School of Industrial and Systems Engineering, and an Associate Director of the Machine Learning Center. She received her Ph.D. in Electrical Engineering (minor in Mathematics) from Stanford University, M.Sc. in Electrical and Computer Engineering from the University of Florida, and B.Sc. in Electrical Engineering and Computer Science from the University of Science and Technology of China (USTC). She was a Research Scientist at Duke University. Her research areas are statistics (in particular sequential analysis and sequential change-point detection), machine learning, and signal processing, providing the theoretical foundation and developing computationally efficient and statistically powerful algorithms. She has worked on such problems in sensor networks, social networks, power systems, crime data analysis, and wireless communications. She received the National Science Foundation (NSF) CAREER Award in 2017. She is currently an Associate Editor for IEEE Transactions on Signal Processing, Sequential Analysis: Design Methods and Applications, and INFORMS Journal on Data Science, and serves on the Editorial Board of  Journal of Machine Learning Research.

\vspace{0.2in}
\noindent
Rui Zhang is a Ph.D. student of H. Milton Stewart School of Industrial and Systems Engineering, Georgia Institute of Technology. Rui Zhang received a B.S. (2015) in statistics from Sun Yat-sen University and a M.S. (2017) in statistics from University of Michigan. He is interested in change-point detection, machine learning and statistics.

\section*{Acknowledgement}

The work is partially funded by an NSF CAREER Award CCF-1650913, NSF CMMI-2015787, DMS-1938106, DMS-1830210.

\bibliographystyle{plain}
\bibliography{References}

\end{document}